\documentclass[
conference]{IEEEtran}
\IEEEoverridecommandlockouts
\usepackage{cite}
\usepackage{amsmath,amssymb,amsfonts}
\usepackage{algorithmic}
\usepackage{graphicx}
\usepackage{textcomp}
\usepackage{xcolor}
\def\BibTeX{{\rm B\kern-.05em{\sc i\kern-.025em b}\kern-.08em
    T\kern-.1667em\lower.7ex\hbox{E}\kern-.125emX}}

\usepackage{makecell}

\begin{document}

\title{An Efficient Data-Driven Model for Generation Expansion Planning with Short-Term
Operational Constraints\\
\thanks{This work was supported by the Natural Sciences and Engineering Research Council (NSERC) Discovery Grant.}
}

\author{\IEEEauthorblockN{1\textsuperscript{st} Hassan Shavandi}
\IEEEauthorblockA{\textit{Management Science and Engineering} \\
\textit{University of Waterloo}\\
Waterloo, Canada \\
hassan.shavandi@uwaterloo.ca}
\and
\IEEEauthorblockN{2\textsuperscript{nd} Mehrdad Pirnia}
\IEEEauthorblockA{\textit{Management Science and Engineering} \\
\textit{University of Waterloo}\\
Waterloo, Canada \\
mpirnia@uwaterloo.ca}
\and
\IEEEauthorblockN{3\textsuperscript{rd} J. David Fuller}
\IEEEauthorblockA{\textit{Management Science and Engineering} \\
\textit{University of Waterloo}\\
Waterloo, Canada \\
dfuller@uwaterloo.ca}
}

\maketitle

\begin{abstract}
Generation expansion planning models have been useful aids for long-term planning in electricity markets. Recent growth in intermittent renewable generation has increased the need to model non-renewable responses to rapid changes in daily loads, leading research to incorporate unit commitment features into generation expansion models. Such combined models usually contain discrete variables which, along with many details, create long computation times that may take days or weeks. This is impractical for analysts who need to develop, debug, modify and use the model for many alternative runs. We propose a novel data-driven generation capacity planning model, in which generation is aggregated by technology type using clustering-based machine learning techniques. The model includes only the minimal unit commitment constraints necessary to integrate the operational requirements and long-term capacity planning. The variations in generation required to respond to rapid demand fluctuations are estimated from historical data, specifically as maximum rates of change for each generation type. We develop our data-driven model using the data from the province of Ontario, Canada. The developed model is a large-scale linear program that can solve problems in less than one hour on modest computing equipment with credible results.
\end{abstract}

\begin{IEEEkeywords}
Capacity Expansion, Unit Commitment, Data-Driven Model, Ramping Constraints, Uncertainty, Machine Learning
\end{IEEEkeywords}

\section{Introduction}

In the early 1970s, Anderson \cite{b1} examined generation expansion planning (GEP) models used by the electric power industry in several countries and discussed the advantages of linear programming (LP) over nonlinear and dynamic programming, and proposed methods to reduce computation time by half or more. Despite advancements in computer hardware and LP software, challenges with long computation times persist due to increasingly complex models that incorporate more planning periods, regions, stochastic scenarios, and distributed generation capacities.

Some suggest that computing time for a GEP model is not a concern since they are typically solved a few times a year for long-term planning studies; for example, \cite{b2} describes a large stochastic GEP model with 4, 5 or 6 investment periods solved in six hours, three days, and over a week, respectively. However, such lengthy computing times can impede the adoption of a GEP model by planning authorities, who need quick access to solutions for modifications with new data, policy changes, sensitivity analyses, and debugging.

The growth of intermittent renewable generation has led research into using short-term unit commitment (UC) features—such as ramping limits, startup/shutdown costs, and minimum power output times—in long-term models. As noted by Poncelet \cite{b4}, many energy system optimization models (ESOMs) use ad-hoc features to approximate how different generation types respond to demand changes. However, these features are often too approximate, or they require difficult parameter estimation \cite{b5}. Thus, we focus on explicitly integrating UC features into planning models.

The GEP model in \cite{b6} incorporates short-term constraints by using a representative day for each season. Compared to traditional approaches, it shows that renewable energy variability can influence operating decisions, leading to sub-optimal outcomes. From \cite{b9}, clustering algorithms have been used to group days with similar loads to select a representative day from each cluster, improving solution accuracy while preserving computational efficiency. Furthermore, a single target year GEP model in \cite{b10} uses most UC features except minimum up/down times and compares it to a reduced version with only ramping constraints and no binary UC features. The reduced LP solves much faster, cutting computation times from hours to about 20 seconds, with similar investment results. 

This paper presents a multi-year GEP model that incorporates the most essential short-term operational constraints to efficiently compute solutions for practical use. Based on results in \cite{b10}, the model excludes binary variables for daily generator on/off status, omitting UC features like startup/shutdown costs, minimum up/down time limits, and minimum run limits. Instead, it uses approximate capacity and variation limits, defined as ramping limits for groups of similar generators, for daily operations. These limits are estimated from historical decisions of the Independent System Operator (ISO). The proposed model tracks aggregated groups (e.g., nuclear, gas, hydro) rather than individual units, using clustering-based machine learning techniques, to eliminate the need for binary variables. Each group's hourly power output variation is limited by maximum rates derived from actual data. With no discrete variables, and having all linear constraints and objective function, the model is an LP featuring multiple representative days to account for seasonal and daily demand variations and extreme scenarios. 

\section{The Integrated Short-Term Operations with Generation Capacity Expansion Problem}

This section defines the optimization framework of the proposed Integrated Short-Term Operations with Generation Capacity Expansion Problem (ISO-GEP) model and presents the mathematical expressions for the constraints and objective function, as shown in \eqref{deqn_ex2} to \eqref{deqn_ex12}. The parameters are indexed by generation technology type \(k\) (e.g., nuclear, gas), years \(t\) (1 to 20), seasons \(ss\), hours \(h\) of each representative day, types of representative days \(i\) (e.g., weekday, weekend), and scenarios \(s\) (e.g., high, medium, or low demand, and initial generator status). It also outlines the methodology used to estimate some parameters. 

\subsection{Data-Driven Methodology for GEP Models}

To better approximate the grouping of generations and the variations among them, we use visualization and clustering techniques to aggregate generators by technology and represent electricity demand by weekdays, weekends, and seasons. Predictive models then estimate the parameters of the ISO-GEP model, which include: 

\begin{itemize}
\item Average hourly demand for representative day \(i\) each season \(ss\).
\item Average capacity factor for each technology \(k\) over season \(ss\) and representative day \(i\).
\item Average variation in demand and supply over season \(ss\) and representative day \(i\).
\item Initial generation levels during each representative day \(i\) in different scenarios \(s\).
\end{itemize}

Lastly, prescriptive modelling uses the proposed ISO-GEP model, a data-driven LP model to integrate short-term operational constraints with long-term generation capacity planning to recommend a data-informed generation mix for the future of energy market of study.

\subsection{The Optimization Framework of the ISO-GEP Model}

We distinguish the existing capacity \(XE(\omega) \) in MW, and total new capacities \( XN(\omega) \) in MW, defined by the new capacity of technology \(k\) to be installed as \( x(\omega') \) in MW, considering the depreciation factor, \( d(k) \), as follows:
\begin{equation}
\label{deqn_ex1}
XN(\omega) = \sum_{t' \leq t} x(\omega') \left[1 - d(k)(t - t')\right]
\end{equation}
The depreciation parameters for the new installed generators are estimated from the life span of each generation technology and the straight-line method. The index 
\(\omega\) represents technology \(k\) and year \(t\), while \(\omega'\) represents technology \(k\) and year \(t'\) before year \(t\).

For the objective, we use a cost minimization function, where the investment cost \(N(\omega)\) in \$/MW, annual fixed cost \(F(\omega)\) in \$/MW, and expected variable and reserve costs \(C(\omega)\) in \$/MWh and \(V(\omega)\) in \$/MW are defined. The expected values of hourly variable and reserve costs are found using probabilities \(P(ss, i, s)\) for scenario \(s\), season \(ss\), and representative day \(i\). The number of representative days \(i\) at season \(ss\) is \(T(i, ss)\), and the variation up and down amount is \(r(\phi)\) in MW. The generation amount of each technology is represented by \(g(\phi)\) in MWh, where \(\phi\) represents technology \(k\), year \(t\), season \(ss\), day \(i\), hour \(h\), and scenario \(s\). The index \(\theta\) is used to represent \(k\), \(t\), \(ss\), \(i\), \(s\). The LP model is defined as follows:
\begin{align} \label{deqn_ex2}
    \min_{x,g} \sum_{\omega} \Big[ N(\omega) x(\omega) + \sum_{t' \leq t} F(\omega) x(\omega') \Big] \\
    + \sum_{\theta} P(ss, i, s) \Big\{ 
    \sum_{h} T(i, ss) \Big( C(\omega)g(\phi) \notag \\
    + V(\omega)r(\phi) \Big) \Big\} \notag
\end{align}
\begin{align} 
    g(\phi) – K(k, ss, h)(XE(\omega)  \label{deqn_ex3} + XN(\omega) \leq 0, \\
    \forall k \neq GAS, t, ss, i, h, s \notag \\
    g(\phi) + \alpha_{RES}(t)G(t)D(ss, i, h, s) \label{deqn_ex4} \\
    - K(k, ss, h)(XE(\omega) + XN(\omega) \leq 0, \notag \\
    k = GAS, \forall t, ss, i, h, s \notag
\end{align}
\begin{align} 
    g\left(\phi\right)-g\left(\phi|h-1\right)-r\left(\phi\right)\leq 0, \forall \phi \label{deqn_ex5} \\
    g\left(\phi|h-1\right)-g\left(\phi\right)-r\left(\phi\right)\le0, \forall \phi \label{deqn_ex6}
\end{align}
\begin{align} 
    g(\phi)-g(\phi|h-1) \label{deqn_ex7}\\
    -U(k,ss)K(k,ss,h)(XE\left(\omega\right)+XN\left(\omega\right))\le0, \notag \\
    k\in\{Nuc,Gas,Hydro,Biofuel\}, \forall\ t,ss,i,h,s \notag \\
    g(\phi|h-1)-g(\phi) \label{deqn_ex8} \\
    -D(k,ss)K(k,ss,h)(XE(\omega)+XN(\omega))\le0, \notag \\
    k\in\{Nuc,Gas,Hydro,Biofuel\}, \forall\ t,ss,i,h,s \notag
\end{align}
\begin{align} \label{deqn_ex9}
g(\phi|h=0)=I(\phi|h=0)(XE(\omega)+XN(\omega)), \\
\forall\ k,t,ss,i,s \notag
\end{align}
\begin{align} \label{deqn_ex10}
G(t)D(ss,i,h,s)-\sum_{k}g(\phi)\le0,
\forall t,ss,i,h,s
\end{align}
\begin{align} 
XE\left(\omega\right)+XN\left(\omega\right)- \label{deqn_ex11} \\
\alpha_H(\omega)\sum_{k}{\left[XE\left(\omega\right)+XN\left(\omega\right)\right]\le0, 
\forall \omega} \notag \\
XE\left(\omega\right)+XN\left(\omega\right)- \label{deqn_ex12} \\
\alpha_L(\omega)\sum_{k}{\left[XE\left(\omega\right)+XN\left(\omega\right)\right]\geq0,\forall \omega} \notag
\end{align}
All cost parameters depend on year \(t\), as costs are discounted to allow for real cost increases. 

Constraints \eqref{deqn_ex3} and \eqref{deqn_ex4} limit the generation output of each type of generation \(k\) to the maximum effective capacity of the generators in each year, season, representative day, hour, and scenario, using capability factor \( K(k, ss, h) \). Gas generators are considered separately, since they are used as reserve capacities in our model; thus, part of their capacity should be put aside for providing reserve ancillary services in each year.
From \eqref{deqn_ex4}, \(\alpha_{RES}(t)\) is the fraction of total demand to be considered as reserve capacity among the gas generation units at year \(t\). The daily load curve is \( D(ss, i, h, s) \) for scenario \(s\) during season \(ss\) and representative day \(i\). When multiplied by growth parameter \( G(t) \), we get the forecast of the load curve for year \(t\). We assume values of \( G(t) \) based on a growth rate of 0.3\% per year; this is slightly less than the middle of the range of possible future demands for Ontario, in the IESO report \cite{b12}.

Constraints \eqref{deqn_ex5} and \eqref{deqn_ex6} define the variation variables that indicate the change in generation between consecutive hours in the objective function.

Constraints \eqref{deqn_ex7} and \eqref{deqn_ex8} represent the variation up and down for generation technology type \(k\) in scenario \(s\).
Where \( U(k, ss) \) and \( D(k, ss) \) are the variation up and down limits of generation group \(k\) at season \(ss\). We estimate these from IESO data on power outputs of generation groups, as results show the engineering estimate of single-unit parameters may be misleading when used for our models. To estimate, we calculate the rate of change of generation output for each generation group over each season, using the maximum positive rate for the variation up and the maximum of absolute negative rate for the maximum variation down. These rates are then estimated as fractions of the group's capable capacity. 

Equation \eqref{deqn_ex9} represents the initial status of generation at hour zero. It is fixed based on the scenarios developed for each year, season, and day.

To satisfy demand in each scenario, we define constraint \eqref{deqn_ex10} where the total supply and demand is balanced every hour of each combination of season, representative day, and scenario. Yearly demand is approximated using growth parameter \( G(t) \).

Since governments maintain different generation target mechanisms, we add constraints \eqref{deqn_ex11} and \eqref{deqn_ex12} on generation related to each technology type. We define parameters \( \alpha_H(\omega) \) and \( \alpha_L(\omega) \) as the maximum and minimum share of each technology group in total generation capacity.

Other constraints are possible, but the model above is enough to see whether the use of variation alone can direct investment in generation that can respond to rapid changes in demand, while maintaining practical computation time.  

\section{Descriptive and Predictive Analysis: Case of Ontario}

The Independent Electricity System Operator (IESO) supply data is organized into six types \(k\), i.e., nuclear, gas, hydro, wind, solar, and biofuel. Public hourly demand data from IESO \cite{b13} from 2003 to 2024 is used to estimate model parameters. Our goal is to cluster historical data to identify patterns, reduce computation time, and ensure realistic demand-supply variability that necessitates investment in ramping capability and generation capacity. 

\subsection{Clustering of Power Demand and Supply Data}\label{3A}

Earlier studies distinguish weekends from weekdays \cite{b2}, \cite{b10}, \cite{b7}, \cite{b8}, \cite{b11}. We analyzed whether Saturday and Sunday power demands can be grouped by calculating the average hourly demand, and results showed similar demands for both days. Thus, we combine them as a single representative weekend day. Results showed significant differences between the average demand for weekdays and weekends. In fact, power demand trends on weekday holidays are the same as weekends; thus, we consider the weekday holidays in the category of weekends. Additionally, seasonal patterns are observed in the average hourly power demand of the days in a year. This measure roughly indicates the seasons ss, defined in Table I. 

\begin{table}[htbp]
\caption{Definition of Seasons\label{tab:table1}}
\centering
\begin{tabular}{l l}
\hline
Seasons & Months\\
\hline
Winter & 2nd half of Nov, Dec, Jan, Feb, Mar\\
Spring & Apr, May, and 1st half of Jun\\
Summer & 2nd half of Jun, July, Aug, and 1st half of Sep\\
Fall & 2nd half of Sep, Oct, and 1st half of Nov\\
\hline
\end{tabular}
\end{table}

Thus, we propose a clustering design including four seasons, and a representative day each for weekdays and weekends, for eight clusters altogether.

\subsection{Prediction of Electricity Demand}
To consider the stochastic variation of hourly power demand on representative day \(i\) and season \(ss\), we use the scenario approach. Three different scenarios \(s\) (low, most likely, high) are created for each \(ss-i\) combination by depicting the demand data as a three-bin histogram. The median of each bin defines a scenario, and its probability is calculated using the frequencies.  

In the model of Section II, the initial “hour zero” power output, defined to be midnight, of generators influences the power produced in coming hours. Since power generation at midnight is similar for weekdays and weekends, the initial status is the same. For each season \(s\) and generator \(k\), initial output is approximated as a fraction of installed capacity and represented as a 2-bin histogram. Only one bin is adequate for solar and for biofuel: leaving two bin histograms for nuclear, gas, hydro, and wind. 

The definition of scenarios combines demand variations and initial generation levels. With four generation types having two hour-zero states (high or low) and three daily load types (high, medium, low), the total number of scenarios, \(s\), is \(2^4 \times 3 = 48\) for each combination of season, \(ss\), and representative day, \(i\).

\section{Results}

This section applies the model of Section II with realistic parameters for Ontario’s electricity market, to show that the model produces credible results in practical time, and that the variation constraints \eqref{deqn_ex7} and \eqref{deqn_ex8} effectively direct investment to generation that can respond to rapidly changing demand.

\begin{table}[htbp]
\caption{Approximate Costs Based on 2016 \label{tab:table2}}
\centering
\begin{tabular}{l l l l l l}
\hline
& \makecell{Life \\ Span \\ (Years)} & \makecell{Invest \\ (\$/MW) \\for Lifetime} & \makecell{Invest \\ (\$/MW) \\ for 20 \\ Years} & \makecell{Fixed \\ O\&M \\ (\$/MW-\\yr)} & \makecell{Variable \\ Cost \\ (\$/MWh)}\\
\hline
NUC & 60 & 7,371,800 & 2,457,267 & 124,347 & 2.852 \\
Hydro & 80 & 7,101,108 & 1,775,277 & 155,000 & 0 \\
GAS & 30 & 1,290,344 & 860,229 & 13,640 & 6.3364 \\
Wind & 30 & 2,327,480 & 1,551,653 & 49,228 & 0 \\
Solar & 30 & 3,244,253 & 2,162,836 & 29,016 & 0 \\
Biofuel & 25 & 6,181,400 & 4,945,120 & 52,204 & 5.208 \\
\hline
\end{tabular}
\end{table}

We consider a planning horizon of 20 years. The life spans of new units in Table II \cite{b19}, \cite{b20} are used to find the straight-line depreciation rate, with variation costs assumed to be 10\% of variable costs. The “Invest for 20 years” column shows how investment costs are allocated based on the portion of the new capacity’s life within the model’s 20-year planning horizon.

\subsection{Performance of the Proposed Model}
The proposed model is a large-scale LP with 829,560 variables and 5,012,400 constraints, solved efficiently by CPLEX in GAMS in 3035 seconds (51 minutes) on a 9-year-old server. A faster computer, and less scenarios \cite{b12}, \cite{b22} could cut the execution time considerably, making the model more useful for policy analysis. Fig. 1 summarizes the main results.

\begin{figure}
    \centering
    \includegraphics[width=0.8\linewidth]{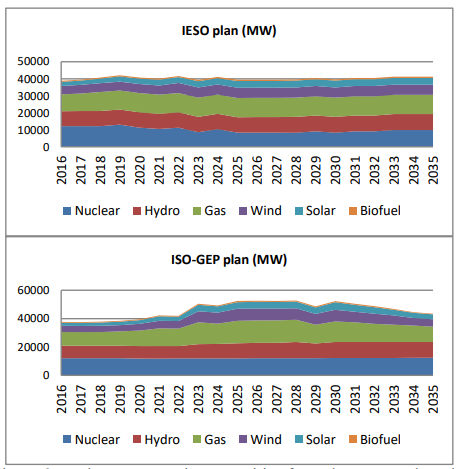}
    \caption{Comparative results on generation capacities between IESO plan and 
proposed model}
    \label{fig:enter-label}
\end{figure}

The difference in supply capacities between our model and the IESO is due to the inclusion of storage capacities and demand response in the IESO-GEP model, which contributes 1416 MW in supply capacities in 2035. Additional discrepancies arise from the assumptions about future demand forecasting. Despite this, our results remain credible and reasonably align with the long-term capacity plan of the IESO.

The share constraints \eqref{deqn_ex11} and \eqref{deqn_ex12} significantly impact the model. Without these constraints, the optimal solution involves investments only in new nuclear and gas, and a small amount in biofuel in the last three years. With the assumed data on costs, the hydro, wind and solar investments are too costly. In another experiment, most share constraints \eqref{deqn_ex11} and \eqref{deqn_ex12} were dropped but kept the lower limit constraints \eqref{deqn_ex12} for wind, solar and biofuel generation. The results in new capacity were closer to the results of the full ISO-GEP model, but again no investment in new hydro was made. 

\subsection{The Impact of Variation Constraints and Costs on Capacity Investments}

To assess the impact of short-term variation costs and constraints \eqref{deqn_ex7} and \eqref{deqn_ex8} on long-term capacity planning, we test when the constraints are binding and evaluate daily operations of investments from models with no variation constraints or variation costs. For the latter, we modify the ISO-GEP model by removing variation costs or constraints (the “no-variation” versions), and then re-run with investments fixed at the no-variation solution, but with variation costs and constraints imposed. We also construct a conventional generation expansion model, with demand aggregated into base, medium, and peak blocks in each year-season-scenario combination; then we run the ISO-GEP model with variation constraints and costs, and with investments fixed at the conventional model’s levels.

\subsection{Binding Variation Constraints in the ISO-GEP Model}

To assess variation constraints in the ISO-GEP model, we examine dual variables of the variation constraints for nonzero values. The variation-up constraint is mostly binding during hours 6-8, and 18, primarily in high-demand scenarios in most seasons and years, due to rising power demand. The variation-down constraint is primarily binding in early morning hours, when demand declines, mostly in low and medium-demand scenarios for most seasons and years. Removing non-binding variation constraints only slightly reduced computing time, from 51 to 50 minutes.

\subsection{Testing Investments Determined by ISO-GEP Without Variation Constraints or Costs}

Variation constraints are further investigated by solving the ISO-GEP model in three cases: no variation constraints or variation costs, no variation constraints with variation costs in the objective, and variation constraints without variation costs. 

Without variation constraints, investment in capacity is much lower, particularly in nuclear. This is because the variation constraints \eqref{deqn_ex7} and \eqref{deqn_ex8} relate the variation limit to total capacity, necessitating investment in new capacity to ensure the supply system can respond to demand variations. 

The two cases—no variation constraints and no variation cost, and no variation constraints with variation cost in the objective—provide the same solution. In the full ISO-GEP model, we fix the capacities and run the model to determine feasibility and cost. The solutions of both cases are infeasible in the fixed-capacities model. In contrast, the model with variation constraints but no variation cost is feasible in the fixed capacities. The difference between this model and the full model is small, but the full model performs slightly better in terms of total cost, as expected.

\subsection{Testing Investments Determined by a Conventional Generation Expansion Model}

The conventional generation expansion model reorders annual demands into a load duration curve, sorting hourly demands in descending order and approximating them by several blocks. Due to the reordering and approximation, it cannot represent hour-by-hour variations in demand and supply, so variation constraints cannot be represented.

We extended the conventional generation expansion model by defining separate three-block approximations for each year-season-scenario combination, as in the ISO-GEP model, and by including the reserve requirement and share constraints. As in the previous subsection, we determine investments with the extended conventional model, then run the full ISO-GEP model with the investments fixed at the extended conventional model’s values. Once again, the result is infeasibility. 

\subsection{Summary of Impacts of Variation Constraints and Costs on Generation Investment}

We conclude that variation constraints are essential for investment decisions in the ISO-GEP model, as they are often binding in the optimal solution, and investments without them are infeasible. Similarly, investment decisions from the extended conventional generation expansion model are also infeasible. In contrast, variation costs have a minor effect on investment decisions.

\section{Summary and Directions for Future Research}

Conventional generation capacity planning models simplify computation by grouping similar generation units in broad categories (e.g., nuclear, gas) and ignoring short-term operational limits and costs by aggregating time into year- or season-long basic units. This overlooks the need for a generation mix that provides capacity to meet peak demands, reserve requirements, and handles rapid daily demand fluctuations.

This paper presents a model that retains the aggregation of generation units but uses 24-hour units for each representative day, each year. The model has no binary variables for investment or generator on/off status and represents limitations on the system to respond to hourly changes in demands by variation up and down constraints for aggregated generation groups. The resulting LP can be solved efficiently for policy analysis, and credible parameter estimates are achieved by analyzing historical operational data for Ontario, Canada. Variation limits for aggregated generation groups are estimated with the assumption that the total limit is proportional to the total installed capacity. 

Future research could expand the model to include storage investment for daily arbitrage and reserves to address uncertainties in renewable supply. Additional features, such as emission limits or taxes, can also be integrated. Further reductions in computation time could be achieved by refining the selection of representative days and scenarios \cite{b12}, \cite{b22}.

\section*{Acknowledgment}

We are grateful for the financial support for this research from NSERC of Canada, through grants to D. Fuller and M. Pirnia, and contributions from Sheng Fang, graduate Masters student at the University fo Waterloo.



\end{document}